\documentstyle [12pt]{article}
\begin{document}

\font\tfont= cmbx10 scaled \magstep3
\font\afont= cmcsc10 scaled \magstep2
\centerline{\tfont On Uniform Homeomorphisms of the}
\smallskip
\centerline{\tfont Unit Spheres of Certain Banach Lattices}
\bigskip
\centerline{\afont F. Chaatit}
\bigskip
%\title{On Uniform Homeomorphisms Of The Unit Spheres of Certain Banach
%       Lattices.}
%\author{F. Chaatit}
%\maketitle
\begin{abstract}

{We prove that if $X$ is an infinite dimensional Banach
 lattice with a weak unit then there exists
a probability space $(\Omega, \Sigma, \mu)$ so that the unit sphere
$S(L_1(\Omega, \Sigma, \mu)$ is uniformly homeomorphic to the unit sphere $S(X)$
if and only if $X$ does not contain $l_{\infty}^n$'s uniformly.}
\end{abstract}

\newcommand{\del}{\mbox{$\delta$}}
\newcommand{\lftyn}{\mbox{${l_{\infty}^{n}}$}}
\newcommand{\norm}[1]{\mbox{$\|#1\|$}}             

\newcommand{\reals}{\mbox{\bf R}}
\newcommand{\integers}{\mbox{\bf N}}
\newcommand{\smint}{\mbox{\scriptsize {\bf N}}}
\newcommand{\complex}{\mbox{\bf C}}
\newcommand{\ep}{\mbox{$\varepsilon$}}
\newcommand{\pf}{\noindent{\bf Proof\,\ }}
\newcommand{\ra}{\mbox{$\longrightarrow$}}
\newcommand{\bdfn}{\noindent{\bf Definition\,} \begin{em}}
\newcommand{\edfn}{\end{em}}
\newcommand{\ntn}{\noindent{\bf Notation\,\,}}
\newcommand{\rmk}{\noindent{\bf Remark\,\,}}
\newcommand{\osm}{\mbox{$(\Omega,\Sigma,\mu)$}}
\newcommand{\om}{\mbox{$\omega$}}
\def\baselinestretch{2}

\newtheorem{thm}{Theorem 2}[chapter]
\newtheorem{lem}[thm]{Lemma 2}
\newtheorem{cor}[thm]{Corollary 2}
\newtheorem{prop}[thm]{Proposition 2}         

\newtheorem{ex}{Exercise}                            

\newenvironment{proof}{\medskip \par \noindent {\bf Proof}\ }{\hfill $\Box$
                       \medskip \par}

\setcounter{page}{1}

\section{Introduction}

Recently E.Odell and Th.Schlumprecht \cite{bib:O.S}
   proved that if $X$ is an
infinite dimensional Banach space with an unconditional
basis then the unit sphere of $X$ and the unit sphere of $l_1$ are
uniformly homeomorphic if and only if $X$ does not contain
$l_{\infty}^{n}$ uniformly in $n$. We extend this result to the setting of
Banach lattices. In Theorem 2.1 we obtain that if $X$ is a Banach lattice
with a weak unit then there exists a probability space $(\Omega, \Sigma, \mu)$
so that the unit sphere $S(L_{1}(\Omega, \Sigma, \mu)$ is
uniformly homeomorphic
to the unit sphere $S(X)$ if and only if $X$ does not contain $l_{\infty}^n$
uniformly in $n$.
A consequence of this -Corollary 2.11- is that if $X$ is a separable
infinite dimensional Banach lattice then $S(X)$ and
$S(l_1)$ are uniformly homeomorphic if and only if
$X$ does not contain $l_{\infty}^n$ uniformly in n.
Quantitative versions of this corollary are given in Theorem 2.2
and  Theorem 2.3. A continuous function $f:[0, \infty) \rightarrow [0,\infty)$
with$f(0) = 0$ is a {\em modulus of continuity} for
a function between two metric spaces      $F:(A, d_1) \rightarrow (B,d_2)$
if $d_{2}(F(a_{1}), F(a_{2})) \leq f(d_{1}(a_{1},a_{2}))$ whenever
$a_{1}, a_{2} \in A$.
Theorem 2.2 says that if $X$ and $Y$ are separable
infinite dimensional Banach lattices with $M_{q}(X) < \infty$
and $M_{q'}(Y) < \infty$
for some $q, q' < \infty$ then there exists a uniform homeomorphism
$F:S(X) \rightarrow S(Y)$ such that $F$ and $F^{-1}$ have modulus of
continuity  $f$ where $f$ depends solely on $q,q', M_{q}(X)$
and $M_{q'}(Y)$.
Here $M_{q}(X)$ is the $q$-concavity constant of $X$ and will be defined below.

Central in defining these homeomorphisms is the entropy map, considered
in \cite{bib:G} and \cite{bib:O.S}. We refer the reader to \cite{bib:B}
and its references
for a  survey of some results concerning
  uniform homeomorphisms between Banach spaces.
In particular it is interesting to note Enflo's result that $\l_1$ and $L_1$
are not uniformly homeomorphic \cite{bib:B} while their unit spheres are.
Also we refer to \cite{bib:L.T} for  facts related to
the theory of Banach lattices.
\vspace{.4cm}

{\bf Aknowledgements:}
 I would like to express my gratitude to
Professors E.Odell and Th.Schlumprecht for proposing this work to
me and for providing me with valuable suggestions and references.
Thanks are also due to Professor V.Mascioni for simplifying the proof of
lemma 2.6.
\vspace{.4cm}

After this work was done, Nigel Kalton discovered a proof of our main
result using complex interpolation theory \cite{bib:K}.
\vspace{.4cm}

\ntn

Let us start by recalling some definitions and well known
facts. A non negative element $e$ of a Banach lattice $X$ is a
{\em weak unit} if $e \wedge x = 0$ for $x \in X$ implies that $x= 0$.
Every separable Banach lattice has a weak unit \cite[p 9]{bib:L.T}.
A Banach lattice is order continuous if and only if every increasing,
order bounded sequence is convergent.
By a general representation theorem (see \cite[p 25]{bib:L.T})
any order continuous Banach lattice with a weak unit can be represented
as a Banach lattice of functions. More precisely:
\begin{enumerate}
 \item there exist a probability space $\osm$
       and an ideal $\widetilde{X}$ of $L_{1}(\Omega,\Sigma ,\mu)$,
       along with a lattice norm $\norm{\cdot}_{\widetilde{X}}$ 

       on $\widetilde{X}$      

        so that $X$ is order isometric
       to $(\widetilde{X}, \norm{\cdot}_{\widetilde{X}})$
 \item $\widetilde{X}$ is dense in $L_{1}(\Omega, \Sigma, \mu)$ and
       $L_{\infty}(\Omega, \Sigma, \mu)$ is dense in $\widetilde{X}$.
 \item ${\norm{f}}_{1} \leq \norm{f}_{\widetilde{X}} \leq 2  \norm{f}_{\infty}$
       for all $f \in L_{\infty}(\Omega, \Sigma, \mu)$.
\end{enumerate}

Moreover ${\widetilde{X}}^* = \{g:\Omega \ra \reals  :
\norm{g}_{\widetilde{X}^*}
 < \infty \}$ is isometric to $X^*$, where
$$\norm{g}_{\widetilde{X}^*} = \sup \{ \int fg d\mu; \norm{f}_{\widetilde{X}}
\leq 1\}$$
and if $g \in \widetilde{X}^*$ and $f \in \widetilde{X}$ then
$$ g(f) = \int fg d\mu.$$

If $X$ is a Banach lattice which is not order continuous then $X$
contains $c_0$
(\cite[pages 6--7]{bib:L.T}).

A Banach lattice $X$ is {\em $q$-concave} if there exists a constant
$M_q < \infty$ such that

$$ (\sum_{i = 1}^{n}{\norm{x_i}}^{q})^{\frac{1}{q}}
\leq M_q \|{(\sum_{i = 1}^{n} |x_i|^{q})^{\frac{1}{q}}}\| \;\;   (\star)$$
(resp.{\em $p$-convex } if there exists $M^p < \infty$
so that
$$ \|{(\sum_{i = 1}^{n}|x_i|^{p})^{\frac{1}{p}}}\| \leq M^{p}
(\sum_{i = 1}^{n} \norm{x_i}^{p})^{\frac{1}{p}} \;\; (\star \star))$$
for all $n \in \integers$ and $x_i \in X$ , $1 \leq i \leq n$.

$M_{q}(X)$  is the smallest constant satisfying $(\star)$ and
$M^{p}(X) $ is the smallest constant that satisfies $(\star \star)$.

Given a Banach lattice of functions $X$, the $p$-convexification
$X^{(p)}$ of $X$ is given by
$$ X^{(p)} = \{ f: \Omega \ra \reals : |f|^p \in X \}$$
with $$ ||| f ||| = {\norm{|f|^p}}^{\frac{1}{p}}.$$

The space $X^{(p)}$ is a Banach lattice with $M^p(X^{(p)}) = 1$
(\cite[p 53]{bib:L.T}).

We will also need the following result. If $X$ is $r$-convex and $s$-concave,
for $1 \leq r,s \leq \infty$ then $X^{(p)}$ is $pr$-convex and $ps$-concave
with
$$ M^{pr} (X^{(p)}) \leq ( M^r (X))^{\frac{1}{p}}$$
and
$$ M_{ps} (X^{(p)}) \leq (M_s (X))^{\frac{1}{p}}.$$

(See \cite[p 54]{bib:L.T}).

We will use standard Banach space notations,
$Ba X = \{ x \in X: \norm{x} \leq 1\}$ will denote the unit ball of $X$
and $S(X) = \{ x \in X: \norm{x} = 1 \}$ the unit sphere of $X$.
If $h$ is a real function on $\Omega$, then
$\mbox{supp}  h = \{\omega \in \Omega: h(\omega) \neq 0 \}$ is the support
of $h$.
If $B \subset \Omega$, then
$ Bh(\omega) = h(\omega) \chi_B (\omega)$ where $\chi_B$
is the indicator function of $B$.

\section
{The main result:}

We now state the main result of this work

\begin{thm}
Let $X$ be an infinite dimensional Banach lattice with
a weak unit. Then there exists a probability space $\osm$
so that $S(L_{1}\osm)$ is uniformly homeomorphic to $S(X)$
if and only if $X$ does not contain $\lftyn$ uniformly in $n$.
\end{thm}

Our proof of Theorem 2.1 will yield two quantitative results:

\begin{thm}
If $X$ and $Y$ are separable infinite dimensional Banach
lattices with $M_q (X) < \infty$ and $M_{q'} (Y) < \infty$
for some $q, q' < \infty$ then there exists a uniform
homeomorphism
$F: S(X) \ra S(Y)$ such that $F$ and $F^{-1}$ have modulus
of continuity $\alpha$ where $\alpha$ depends solely on
$q, q', M_q (X)\;\mbox{and}\; M_{q'} (Y).$
\end{thm}

\begin{thm}
If $X$ and $Y$ are both uniformly convex and uniformly smooth
separable infinite dimensional Banach lattices then there exists
a uniform homeomorphism $F:S(X) \ra S(Y)$ such that $F$ has
modulus of continuity $f$ where $f$ depends solely on the modulus of
uniform convexity of $Y$ and the modulus of uniform smoothness of $X,$
and $F^{-1}$ has a modulus of continuity $g$ depending solely on the modulus
of uniform smoothness of $Y$ and the modulus of uniform convexity of $X.$
\end{thm}

The proofs  will involve a sequence of steps similar
to those in \cite{bib:O.S}.
We begin with a simple extension of Proposition 2.8 of \cite{bib:O.S}.
Recall that $X^{(p)}$ is the $p$-convexification of $X$.

\begin{prop}
Let $X$ be a Banach lattice of functions on a set $\Omega$ and
let $1 < p < \infty$.
Then the map $$G_p : S(X^{(p)}) \ra S(X)$$ given by
$G_{p} (f) = |f|^{p} \mbox{sign} f$ is a uniform homeomorphism. Furthermore
the moduli of continuity of $G_p$ and $(G_p)^{-1}$ are functions solely of $p$.
\end{prop}

%that was prop2.4

\pf
Clearly $G_p$ maps $S(X^{(p)})$ one-to-one onto
$S(X).$ Let $f$ and $g$ be in $S(X^{(p)})$ with
$1 > \delta = \norm{f-g}_{X^{(p)}} = \norm{|f-g|^p}_X^{\frac{1}{p}}.$

{\noindent As} in \cite{bib:O.S} we shall show that there exist two functions
$H$ and $F$ such that
$$H(\delta) \leq \norm{G_p (f) - G_p (g)} \leq F(\delta) $$

{\noindent where}
 $F(\delta) = 2 (1 - (1 - \del^{\frac{1}{p}})^{p}) + \del^{p-1} +
\del^{p}$ and $H(\del) = \frac{1}{2^{p-1}}  \del^{p}$.
The proposition then follows.

{\noindent Let}  $$\Omega _{+} = \{ \om \in \Omega : \mbox{sign} f(\om) = \mbox{sign}
 g(\om) \}$$
and $$\Omega _{-} = \{ \om \in \Omega : \mbox{sign} f(\om) \neq \mbox{sign}
 g(\om) \}.$$
We then have:
\begin{eqnarray*}
  \norm{G_p (f) - G_p (g)} &=&  \|{| |f|^p  \mbox{sign} f
                               - |g|^p \mbox{sign} g |}\| \\
                           &=&  \|{| |f|^p - |g|^p | \chi_{\Omega_{+}}
                                  + (|f|^p + |g|^p ) \chi_{\Omega_{-}}}\|.
\end{eqnarray*}
But $a^p - b^p \geq (a-b)^p$
and $a^p + b^p \geq 2^{1-p} (a+b)^p$ for $a \geq b \geq 0$.

{\noindent Thus,}
\begin{eqnarray*}
  \left\|G_p (f) - G_p (g)\right\| & \geq &
   \left\|{\left||f|-|g|\right|^p \chi_{\Omega_{+}} 
                                        + \frac{1}{2^{p-1}} ( |f| + |g| )^p
                                        \chi_{\Omega_{-}}}\right\| \\
                           & \geq &
   \left\|{\frac{1}{2^{p-1}}\left||f|-|g|\right|^p \chi_{\Omega_{+}}
                                    +\frac{1}{2^{p-1}} (|f| + |g|)^p
                                        \chi_{\Omega_{-}}}\right\| \\
                           &  =   & 2^{1-p} \|{|f-g|^p}\| \\
                           & = & 2^{1-p} \|{f-g}\|^p_{X^{(p)}} \\
\end{eqnarray*}

{\noindent So} we
obtain $H(\del) = \frac{1}{2^{p-1}} \del^p$ as a lower estimate.
For the upper estimate we have:
\begin{eqnarray*}
  \norm{G_p (f) - G_p (g)} & = & \norm{\left| |f|^p - |g|^p \right|
                                   \chi_{\Omega_{+}}
                                  + (|f|^p + |g|^p) \chi_{\Omega_{-}}} \\
                           & \leq & \norm{| |f|^p - |g|^p |\chi_{\Omega_{+}}}
                            + \norm{ (|f|^p + |g|^p) \chi_{\Omega_{-}}} \\
\end{eqnarray*}

{\noindent First} we note that since $$(|f|^p + |g|^p) \chi_{\Omega_{-}} \leq
 (|f| + |g|)^{p} \chi_{{\Omega}_{-}} \leq \left|f - g\right|^p \chi_\Omega,$$
we get
$$\left\| (|f|^p + |g|^p) \chi_{\Omega_{-}}\right\| \leq
\|f - g\|^p_{X^{(p)}} = \delta^p. $$

{\noindent Next}
 we estimate $\| \left| |f|^p - |g|^p \right| \chi_{\Omega_{+}}\|.$
For this purpose we split $\Omega_{+}$ into $\Omega^1_{+}$ and $\Omega^{2}_{+}$
where
$$\Omega^{1}_{+} = \left\{ \om \in \Omega_{+} : |f(\om) \leq q |g(\om)|
\; \;   \mbox{or} \; \; |g(\om)| \leq q |f(\om)| \right\}$$
and
$$\Omega^{+}_{2} = \Omega_{+} \sim \Omega^1_{+}$$
and $q = 1 - \delta^{\frac{1}{p}}.$

{\noindent Note} that if $C = ( 1 - q )^{-p}$ then
$$| |f|^p - |g|^p | \chi_{\Omega^{1}_{+}} \leq C |f - g |^p.$$
Indeed,
$$C|f - g|^p - |g|^p + |f|^p \geq C | g - qg|^p - |g|^p = 0$$
in case $ |f| \leq q |g|$ (the proof is similar if $|g| \leq q |f|).$

{\noindent Thus}
\begin{eqnarray*}
 \| \chi_{\Omega^1_+}| |f|^p - |g|^p | \| & \leq & C \| \chi_{\Omega_{+}}
              |f - g|^p \| \\
                                          & \leq & C \| |f - g|^p \| \\
                                          & = & C \| f - g \|^p_{X^{(p)}} \\
                                          & = & C \delta^{p} \\
                                          & = & (1 - q )^{-p} \delta^{p} \\
\end{eqnarray*}

{\noindent And} since $( 1 - q)^{-p} = \delta^{-1},$ we obtain
$$\| \chi_{\Omega^{1}_{+}} | |f|^p - |g|^p | \| \leq \delta^{p-1}.$$
Finally we have on $\Omega^2_{+}:$
\begin{eqnarray*}
  \| | |f|^p - |g|^p| \chi_{\Omega^2_{+}} \|
       & \leq & ( 1 -q^p) \| |f|^p + |g|^p | \| \\
       & \leq & 2 (1 - (1 - \delta^{\frac{1}{p}})^p). \\
\end{eqnarray*}

{\noindent So}
 $$ F(\delta) = 2 (1 - ( 1 - \delta^{\frac{1}{p}})^p) + \delta^{p-1} +
\delta^p$$
and as $p > 1,$ $F(\delta) \ra 0$ when $\delta \ra 0 \bullet$

\vspace{1 cm}
Throughout the rest of the paper, $X$ will denote  a Banach lattice with the representation  
as a lattice of functions on $\osm$ satisfying the
conditions mentionned in the introduction.
The next step in proving Theorem 2.1 will be to produce a uniform homeomorphism
$$F_X : S(L_1\osm) \ra S(X) $$
in the case where our lattice $X$ is uniformly convex and uniformly
smooth.
In order to do this we need first to define the {\em entropy function}
$E(h,f).$

Let $h \in (L_{\infty}(\mu))^{+}$ and define
 $E(h,\cdot) : X \ra \bar{\reals}$  by
$$E(h,f) = \int h \log |f| d\mu$$ for $f \in X,$
(we use the convention that $0 \log 0 \equiv 0$)
and more generally,
$$E(h,f) = E( |h| , |f| )$$ if $h  \in {L_{\infty}(\mu)}.$

The entropy map was considered in \cite{bib:G} and in the sequel
we use arguments of both \cite{bib:O.S} and \cite{bib:G}.

\begin{prop}
Suppose $X$ is uniformly convex. Let $h \in (L_{\infty}(\mu))^{+} $
and set $$\lambda \equiv \sup_{f \in Ba X}\int h \log |f|d\mu.$$
Then $-\log 2  \leq \lambda \leq {\|h\|}_{\infty}$ and
if $h \neq 0$ there exists a unique  $f \in S(X)^{+}$
so that $\lambda = E(h,f).$
Moreover  $supp f = supp h.$
\end{prop}
\pf
First we note that $\lambda \leq {\|h\|}_{\infty}.$
To see this it suffices to observe
that
\begin{eqnarray*}
 \lambda & = & \sup_{f \in Ba X^+}\int h \log | f| d\mu \\
         & \leq & \sup_{f \in Ba X^+}\int h |f| d\mu \\
         & \leq & \sup_{f \in Ba X^+} \|h\|_\infty \|f\|_{L_1} \\
         & \leq & \sup_{f \in Ba X^+} \|h\|_\infty \|f\|_X \\
         & \leq & \|h\|_\infty. \\
\end{eqnarray*}
Also $\lambda \geq - \log 2$  since $\chi_{\Omega}/2 \in Ba(X)^{+}.$
Next let $(f_n) \subseteq (Ba X)^{+}$ be such that
$E(h,f_n) \geq \lambda - 2^{-n}.$
Since $X$ is uniformly convex, by passing to a subsequence,
we can suppose that $f_n$ converges weakly to $f \in (Ba X)^{+}.$
Let $(u_n)$ be a sequence of ``far-out'' convex combinations of $f_n,$
such that $(u_n)$ converges to $f$ in norm, thus
$u_n = \sum_{i=p_n +1}^{p_{n+1}}c_{i}f_{i}$ where $p_1 < p_2 <... < p_n <....
c_{i} \geq 0, \sum_{i=p_n +1}^{p_{n+1}}c_i = 1$
and  $\|u_n - f\|_X \ra 0$ as
$n \ra \infty.$

We next note that if  $(g_i)_{i=1}^{n} \subseteq Ba X,$
and $(d_{i})_{i=1}^{n} \subseteq (\reals)^{+}$ with $\sum_{i=1}^{n}d_i = 1$
then$$E\left( h , \sum_{i=1}^{n} d_{i}g_{i} \right) \geq
\sum_{i=1}^{n} d_{i}E(h,g_{i}).$$
Moreover if $B = \mbox{supp} h$ and $Bg_{i} \neq Bg_{j}$ for some $i,j$ then
$$E\left( h , \sum_{i=1}^{n} d_{i}g_{i} \right) > \sum_{i=1}^{n}d_{i}E(h,g_i)$$
This follows from the strict concavity of the logarithm function.

{\noindent Therefore} $$\lim_{n \rightarrow \infty} E ( h , u_n ) = \lambda.$$
\vspace{.4cm}

{\bf Claim:} $E ( h , f ) = \lambda$
\vspace{.5cm}

{\noindent Note} that
  $$\|u_n - f \|_{L_{1}(\mu)} \leq \|u_n - f \|_{X} \rightarrow 0$$
and so in order to prove the Claim, it suffices to prove the following lemma:
\begin{lem}
Let $\lambda \in \reals, h \in L_1^{+}(\mu), 
(u_n) \subseteq L_{1}^{+}(\mu)$
and suppose $u_n \ra f$  in $L_{1}(\mu).$
Then $$\int h \log u_n d\mu \ra \lambda \; \;
\mbox{implies} \;\;
\int h \log f d\mu \geq \lambda.$$
\end{lem}
\vspace{.5cm}

\pf
By passing to a subsequence we may assume that $u_n \rightarrow f$ a.e. Thus
$(\log u_{n})^{-} \rightarrow (\log f)^{-}$ a.e. and so
$$\int h(\log f)^{-}d\mu \leq \liminf_{n \rightarrow \infty}
\int h(\log u_n)^{-}d\mu.$$
by Fatou's lemma.
Therefore
$$(\star)\;\; \limsup_{n \rightarrow \infty}\int -h(\log u_n)^{-}d\mu  \leq \int-h (\log f)^{-} d\mu.$$
On the other hand, one has also the inequality:
$$(\star \star)\;\; \limsup_{n \rightarrow \infty}\int h (\log u_n)^{+} d\mu
\leq \int h (\log f)^{+}d\mu.$$
Indeed, fix $\ep > 0.$
Since $0 \leq (\log u_n)^{+} \leq u_n,$ and $(u_n)$ is uniformly integrable,
there exists $\delta > 0$
so that $\mu (A) < \delta$ implies
$$\mbox{for all}\; \; n, \int_{A}(\log u_n)^{+}d\mu < \ep \; \;  \mbox{and} \; \;
\int_{A}(\log f)^{+} d\mu < \ep.$$
($(\log f)^{+}$ is integrable since $0 \leq (\log f)^{+} \leq f.$)
Now $h(\log u_n)^{+} \ra h(\log f)^{+}$ a.e:
So by Egoroff's theorem, there exists  a set $C$ with $\mu (C) < \delta$
such that $$h(\log u_n)^{+} \ra h(\log f)^{+}$$ uniformly
except perhaps on $C.$
More exactly, for $\ep > 0,$ there exist $n(\ep) \in \integers$ and a set $C$
with $\mu (C) < \delta$ such that for any
$n \geq n(\ep)$
we have $$\sup_{\om \in C^{c}}|h(\log u_n)^{+} - h(\log f)^{+}| < \ep.$$
Thus
\begin{eqnarray*}
 \int h(\log u_n)^{+} d\mu & \leq & \int |h(\log u_n)^{+}-h(\log f)^{+}|d\mu
                                    + \int h(\log f)^{+} d\mu. \\
                        & = & \int_{C}|h(\log u_n)^{+}-h(\log f)^{+}|d\mu \\
                           &   & +\int_{C^{c}}|h(\log u_n)^{+} -
                                h(\log f)^{+}|d\mu
                                  +\int h(\log f)^{+}d\mu.  \\
                           & < & 2 \ep + \ep + \int h(\log f)^{+} d\mu \\
\end{eqnarray*}
So $$\limsup_{n \rightarrow \infty}\int h(\log u_n)^{+} d\mu \leq
\int h (\log f)^{+} d\mu.$$
Now adding $(\star)$ and $(\star \star)$ yields
$$\lambda \leq \int h \log f d\mu,$$ which proves Lemma 2.6 $\bullet$

\vspace{1cm}
Note that since $\lambda \geq E(h,f),$ we get $E(h,f) = \lambda,$
proving the Claim.
Now we prove that $f$ is unique.
Indeed, let $f \neq g$ with $E(h,f) = E(h,g)$ and $\|f\| = \|g\| = 1.$
Thus by uniform convexity$\left\| \frac{f+g}{2}\right\| < 1$
and so $\frac{f+g}{2}$ cannot maximize the entropy,
while clearly $\mbox{supp} h \subseteq \mbox{supp} g$ a.e
and so
$$\lambda =\frac{1}{2}\left(E(h,f) + E(h,g)\right)
\leq E(h,\frac{f+g}{2}) < \lambda,$$
a contradiction.

Let now $B = \mbox{supp}h.$ In order to obtain $\mbox{supp}f = B$ a.e
consider first $g=Bf$ in what preceeds to get $f = Bf.$ Then observe
that trivially $\mbox{supp}Bf \subset B$ a.e, while if the previous
inequality was strict, then there exists a set $A \subset B$ with
$\mu (A) > 0$ such that $f_{|A} = 0.$ Thus
$$- \infty = E(h,f) \geq E(h,\chi_{\Omega}/2) = -\log 2;$$
a contradiction. Hence $\mbox{supp}f = \mbox{supp}Bf = B.$
$\bullet$

\vspace{1cm}
Thus under the assumption that $X$ is uniformly convex
we can define
$$F_X : S(L_1(\mu))^{+} \bigcap L_{\infty}(\mu) \ra S(X)^{+}$$
by $F_{X} (h) = f$ where $f \in S(X)^{+}$ is such that
$$E(h,f) = \max_{g \in (Ba X)^{+}}\int h \log |g| d\mu = E_{X}(h)$$
We then define $$F_{X} : S(L_1(\mu)) \bigcap L_{\infty}(\mu) \ra S(X)$$
by $F_{X}(h) = (\mbox{sign}h) F_{X} (|h|).$

We shall show that $F_X$ is uniformly continuous, and thus extends to a
continuous function on $S(L_{1}(\mu)).$
To do so we will need a proposition similar to Proposition 2.3.C
of \cite{bib:O.S}. The proof is nearly the same, adapted to function spaces.

\begin{prop}
Let $h_1, h_2$ be in $S(L_1(\mu))^{+} \bigcap L_{\infty}(\mu)$ with $\|h_1 - h_2\|
_{1} \leq 1.$
Let $x_1 = F_{X} (h_1),$ and $x_2 = F_{X} (h_2).$ Then
$$\left\|\frac{x_1 + x_2}{2}\right\|
\geq 1 -\|h_1 - h_2\|^{\frac{1}{2}}_{1}.$$
\end{prop}
\pf
Let $\left\|\frac{x_1+x_2}{2}\right\| = 1 - 2\ep.$
We need to show that $$2 \ep \leq \|h_1 - h_2\|^{\frac{1}{2}}.$$
We may assume $\ep > 0.$
Define $\widetilde{x_1}  = x_1 + \ep x_2$
and $\widetilde{x_2} = x_2 + \ep x_1.$
Then $$\mbox{supp}\widetilde{x_1} = \mbox{supp}\widetilde{x_2} =
\mbox{supp}h_1 \cup \mbox{supp}h_2 \equiv B,$$
and $$\left\|\frac{\widetilde{x_1} +\widetilde{x_2}}{2}\right\|
\leq \left\|\frac{x_1 + x_2}{2}\right\| + \ep = 1 - \ep$$
With this we  can prove that
$$\ep \leq |\log (1-\ep)| \leq \frac{1}{2}\{E(h_1,\widetilde{x_1}) -
E(h_1, \widetilde{x_2})\}\;\; (\star)$$
Indeed, since $\widetilde{x_1} \geq x_1,$ we clearly have
\begin{eqnarray*}
  E(h_1, \widetilde{x_1})
             & \geq & E(h_1 , x_1) \\
             & \geq & E\left(h_1,\frac{\widetilde{x_1}+
                        \widetilde{x_2}}{2(1-\ep)}\right)\\
\end{eqnarray*}
since $\frac{\widetilde{x_1} + \widetilde{x_2}}{2(1-\ep)} \in Ba X$ and $x_1$ maximizes the  
entropy.
And
\begin{eqnarray*}
  E\left(h_1,\frac{\widetilde{x_1} + \widetilde{x_2}}{2(1-\ep)}\right)
       & = & E\left(h_1, \frac{\widetilde{x_1} + \widetilde{x_2}}{2}\right)
             +|\log (1 - \ep)| \\
       & \geq & \frac{1}{2}E(h_1,\widetilde{x_1}) +
                \frac{1}{2}E(h_1,\widetilde{x_2})
                +|\log(1 - \ep)| \\
\end{eqnarray*}

{\noindent Similarly} we have
$$\ep \leq |\log (1 -\ep)| \leq \frac{1}{2}\{E(h_2,\widetilde{x_2}) -
  E(h_2,\widetilde{x_1})\}\;\; (\star \star).$$
Then by averaging $(\star)$ and $(\star \star)$ we get
$$\ep \leq \frac{1}{4}[E(h_1,\widetilde{x_1}) - E(h_1, \widetilde{x_2})
+E(h_2, \widetilde{x_2}) - E(h_2, \widetilde{x_1})].$$
 So
\begin{eqnarray*}
 \ep & \leq & \frac{1}{4} \int_{B} (h_1 - h_2)(\log \widetilde{x_1} -
                     \log \widetilde{x_2}) d\mu \\
     & \leq & \frac{1}{4} \int_{B} |h_1 - h_2| \left|\log
               \frac{\widetilde{x_1}}{\widetilde{x_2}}\right| d\mu  \\
\end{eqnarray*}

{\noindent But}
$$\left| \log \frac{\widetilde{x_1}}{\widetilde{x_2}}\right| \leq
       \log \frac{1}{\ep} \; \; \mbox{on} \; \; B$$
for$$\frac{\widetilde{x_1}}{\widetilde{x_2}} = \frac{x_1 + \ep x_2}{x_2 + \ep x_1} =
\frac{x_1 + \ep x_2}{\ep (x_1 + {\ep}^{-1}x_2)} \leq \frac{1}{\ep}.$$
and similarly $$\frac{\widetilde{x_2}}{\widetilde{x_1}} \leq \frac{1}{\ep}.$$
Since $\log \frac{1}{\ep} \leq \frac{1}{\ep},$ we finally get
$$\ep \leq \frac{1}{4} \|h_1 - h_2\|_1 \frac{1}{\ep}.$$
Hence $$2 \ep \leq \|h_1 - h_2\|_1^{\frac{1}{2}}
\bullet$$

\begin{prop}
Let $X$ be uniformly convex. Then
$$F_X : S(L_1(\mu)) \bigcap L_{\infty}(\mu) \ra S(X)$$ is
uniformly continuous and hence extends to a uniformly
continuous map $F_X :S(L_1(\mu)) \ra S(X).$
Moreover the modulus of continuity of $F_X$ depends
only on the modulus of uniform convexity of $X.$
\end{prop}

\pf
Recall that $X$ is uniformly convex if and only if
$${\delta}_X(\ep) = \inf \left\{1 - \left\|\frac{x+y}{2}\right\|:
\|x\| = \|y\| = 1, \|x-y\| \geq \ep \right\} > 0.$$
We first observe that $F_{X}: S(L_{1}(\mu))^{+} \ra S(X)$ is
uniformly continuous.

{\noindent Indeed,} by Proposition 2.7,  if $h_1$ and $h_2$ are in
$S(L_1(\mu))^{+} \bigcap L_{\infty}(\mu)$ and
\newline $\|h_1 - h_2\|_1 \leq 1$
then
$$\left\|\frac{F_X(h_1) + F_X(h_2)}{2}\right\| \geq
1 - \|h_1 - h_2\|_1^\frac{1}{2}$$
or
$$1 - \left\|\frac{F_X(h_1) + F_X(h_2)}{2} \right\|
\leq \|h_1 - h_2\|_1^{\frac{1}{2}}.$$
So if $\|F_X(h_1) - F_X(h_2)\| \geq \ep$ then $\|h_1 - h_2\| \geq
(\delta _{X}(\ep))^2.$
Thus there exists $\eta (\ep) = (\delta _{X}(\ep))^2$ so that
$\|h_1 - h_2 \| < \eta (\ep)$ implies $\left\|F_X (h_1) - F_X (h_2)\right\|
\leq \ep.$
Letting $\eta (0) = 0,$ 

the function $\eta$ is continuous and strictly increasing on
$[0,2].$ So $\eta$
has an  inverse $g$ depending only
on the modulus of uniform convexity of $X,$ and
$$\|F_X (h_1) - F_X (h_2) \| \leq g(\|h_1 - h_2 \|).$$

\vspace{1cm}
For the general case let $h_1, h_2$ in $S(L_1(\mu)) \bigcap L_{\infty}(\mu)$
and set
$$x_i = F_X (h_i) = \mbox{sign} h_i \cdot F_X (|h_i|)$$ for $i=1,2.$
Then
$$\|x_1 - x_2\| \leq \|F_X (|h_1|) - F_X (|h_2|)\| +
\|\chi_D (F_X(|h_1|) + F_X(|h_2|))\|$$
where
$$D = \{ \om \in \Omega : \mbox{sign} h_1 (\om) \neq \mbox{sign} h_2 (\om)\}.$$
By what we observed in the beginning of the proof,
\newline $\|F_X (|h_1|) - F_X (|h_2|)\| < g(\ep)$ whenever
$$\| |h_1| - |h_2| \| \leq \|h_1 - h_2 \| < \ep.$$
Our next step is to estimate $\|\chi_D F_{X}(|h_{i}|)\|,\;\;\mbox{for}
 \;\; i=1,2.$
To do so, we note that
$$\|\chi_{D}F_{X}(|h_{1}|)\| = \|DF_{X}(|h_{1}|)\|
\leq \left\| F_{X}(|h_{1}|)-F_{X}
\left(\frac{D^{c}|h_{1}|}{\|D^{c}|h_{1}|\|} \right)\right\|.$$
We are then lead to estimate
\begin{eqnarray*}
 \left\| h_1 - \frac{D^c h_1}{\|D^c h_1\|} \right\|
   & \leq & \left\|D(h_1 - \frac{D^c h_1}{\|D^c h_1\|})\right\|
            +\left\|D^c(h_1 - \frac{D^c h_1}{\|D^c h_1\|})\right\| \\
   & = & \|D h_1 \| + \left\|D^c h_1 - \frac{D^c h_1}{\|D^c h_1\|}\right\|. \\
\end{eqnarray*}

{\noindent We} first get that
$$\|Dh_1\| = \|D|h_1| \| \leq \|D(|h_1| + |h_2|)\| \leq \|h_1 - h_2\| < \ep;$$
and, since $\|h_1\| = \|Dh_1 + D^c h_1 \| = 1$ and $\|D h_1\| < \ep,$
an easy computation yields
$$\left\|D^c h_1 - \frac{D^c h_1}{\|D^c h_1\|}\right\| \leq \|D h_1\| < \ep.$$
So $\left\|h_1 - \frac{D^c h_1}{\|D^c h_1\|}\right\| < 2 \ep$
and thus
\begin{eqnarray*}
 \|DF_X(|h_1|)\| & \leq & \left\|F_X(|h_1|) - F_X
                          \left(
                         \frac{D^c|h_1|}{\|D^c |h_1|\|}
                         \right)
                          \right\| \\
                 & \leq & g(2 \ep). \\
\end{eqnarray*}
Similarly $\|DF_X (|h_2|)\| \leq g(2 \ep).$
Hence $\|F_X(h_1) - F_X (|h_2|)\| \leq g(\ep) + 2 g(2 \ep).$

Therefore $F_X$ extends uniquely to a uniformly continuous map, that
we still denote $F_X,$
from $S(L_1(\mu))$ to $S(X),$ and the modulus of continuity
of $F_X$ depends
only on the modulus of uniform convexity of $X$.
$\bullet$
\begin{prop}
Let $X$ be uniformly convex and uniformly smooth. Then $F_X : S(L_1(\mu))
\ra S(X)$ is a uniform homeomorphism. Moreover
$(F_X)^{-1} : S(X) \ra S(L_1(\mu))$ has modulus of continuity
depending only on the modulus of uniform smoothness of $X.$
Furthermore $(F_X)^{-1}(x) = |x^*| \cdot x $ where $x^* \in S(X^*)$ is
the unique supporting functional of $x.$   

\end{prop}

\pf
Our goal now is to show that the map $F_X$ previously
defined is invertible and that $(F_X)^{-1}$ has the described form
and is uniformly continuous.

\vspace{1cm}
{\noindent \bf Claim 1:}
{\sl Let $h \in S(L_{1}(\mu)) \bigcap L_{\infty}(\mu)$. Then
$g = F_{X}(h)^{-1} \cdot h \in S(X^{*})$ where $\cdot$
denotes the pointwise product.        }

\vspace{1cm}
Note that $\mbox{supp}F_{X}(h) = \mbox{supp}h$ and we define
$F_{X}(h)^{-1} \cdot h$ to be $0$ off the support of $h.$
Assume Claim 1 for the moment.

{\noindent For}
 $x \in S(X)$, define $G(x) = |x^*| \cdot x,$  where $x^*$ is the
unique supporting functional of $x.$
Let $h \in S(L_{1}(\mu)) \bigcap L_{\infty}(\mu).$ Since
$\mbox{sign} F_{X}(h) = \mbox{sign}h$,
$$\int \frac{h}{F_{X}(h)}|F_{X}(h)| d\mu = \int |h| d\mu = 1.$$
Thus from Claim 1 it follows that
$$\frac{h}{F_{X}(h)} = |F_{X}(h)|^{*} = |F_{X}(h)^{*}|.$$
{\noindent Hence} $G(F_X (h)) = |F_{X}(h)|^{*} \cdot F_{X}(h) = h$
for any $h \in S(L_1(\mu))\bigcap L_{\infty}(\mu).$

Furthermore $G$ is uniformly continuous.
Indeed, the support functional $x \mapsto x^{*}$ is uniformly continuous
since $X$ is uniformly smooth,
and since $G(x_i) = |x_i^*| \cdot x_i \;\; i=1,2$ we have
\begin{eqnarray*}
 \|G(x_1)-G(x_2)\| & = & \||x^*_1| \cdot x_1 - |x^*_2| \cdot x_2\| \\
                   & \leq & \||x^*_1| \cdot (x_1 - x_2)\|
                           + \|(|x^*_1| - |x^*_2|) \cdot x_2 \| \\
                   & \leq & \|x_1 -x_2\| + \| x^*_1 - x^*_2\|. \\
\end{eqnarray*}
Thus $G$ is uniformly continuous. Moreover since the modulus of continuity of
$x \mapsto x^*$ depends only on the modulus of uniform smoothness of $X,$
the same is valid for $G.$
Thus $G(F_{X}(h)) = h$ for all $h \in S(L_{1}(\mu)).$

\vspace{1cm}
{\noindent \bf Claim 2:}  {\sl $G$ is one-to-one.}

\vspace{1cm}
It then follows that $G = (F_X)^{-1}.$
We now prove Claim 1

\vspace{1cm}
\noindent{\bf Proof of Claim1:}
We will follow the path of \cite{bib:G}.
Let $h \in S(L_{1}(\mu)) \bigcap L_{\infty}(\mu)$ and suppose $x = F_X (h).$ We can  
assume that
$h \in S(L_1(\mu))^{+} \bigcap L_{\infty}(\mu).$ Then
\newline $\mbox{supp} x = \mbox{supp} h \equiv B
\;\; \mbox{and}\;\; x \in S(X)^{+}.$
Let $k \in X^+$ be arbitrary, then
$$\infty \; >  E(h,x) \geq \int h \log \frac{x+k}{\|x+k\|}d\mu.$$
{\noindent So} writing $x + k = x (1 + \frac{k}{x})$ for $x \in B$ yields
 $$ E(h,x) \geq E(h,x) + \int_B h \log(1 + k x^{-1}) d\mu - \log \|x+k\|.$$
This gives:
\begin{eqnarray*}
 \int_B h\log(1+kx^{-1})d \mu & \leq & \log \|x+k\| \\
                            & \leq & \log (\|x\| + \|k\|) \\
                            & = & \log (1 + \|k\|). \\
\end{eqnarray*}

{\noindent So} $$\int_B h \log (1+kx^{-1}) d\mu \leq \|k\| \;\;\; (\star).$$

{\noindent Thus } on $B, \;\; kx^{-1}$ is finite $\mu$-almost everywhere.
Let $$\sigma_{n} = \{ \om \in B : k(\om)x^{-1}(\om) \leq n   \}$$
and $\chi _n = \chi _{\sigma_n}$ then $\chi_{n} \nearrow \chi_B,$
pointwise $\mu$-a.e;
and since $t \leq \log (1+t) + \frac{1}{2}t^2$ holds
for all $t \geq 0$ we have for $0 < s < \infty$
\begin{eqnarray*}
 s \int_B hx^{-1}k \chi_n d \mu & \leq & \int_B h \log (1+skx^{-1}\chi_n )d\mu
                  + \frac{1}{2}s^2\int_B k^2(x^{-1})^2 \chi_n h d\mu \\
                              & \leq &\int_B h \log (1 + skx^{-1}) d\mu
                   +\frac{1}{2}s^{2} n^{2}\\
                              & \leq & s \|k\| + \frac{1}{2}s^2 n^2 \;\;
                                         \mbox{by} \;\; (\star).
\end{eqnarray*}
Thus dividing by $s$ and letting $s$
go to $0,$ we obtain for all $n \in \integers$
$$\int h x^{-1} k \chi_{n}d\mu \leq \|k\| ;$$
and therefore by the  monotone convergence theorem,
$$\int _Bh x^{-1}k d\mu \leq \|k\|.$$
 Now let $g=hx^{-1}.$ The previous equality yields $\|g\|_{X^*} \leq 1.$
On the other hand
\begin{eqnarray*}
1=\left| \int h d\mu \right| &  =  & \left|\int g \cdot x d\mu \right| \\
                             & \leq & \|x\|_X \|g\|_{X^*} \\
\end{eqnarray*}
{\noindent So} $\|g\|_{X^*} = 1$ which proves Claim 1.
$\bullet$

\vspace{1cm}
{\noindent \bf Proof of Claim 2:}
{\noindent Let} $h = |x_1^*| \cdot x_1 = |x_2^*| \cdot x_2$ be
a member of $S(L_1(\mu))$ with $x_i^*(x_i) = 1,
x_i \in S(X)$ and $x_i^* \in S(X^*)\; \; \mbox{for}\; \; i=1,2.$
We first note that supp$h = $ supp$x_i$ for $i=1,2.$
Indeed $\mbox{supp}h \subset \mbox{supp}x_i$ is clear,
and in case the inclusion is strict
let us consider $B|x_i|$
where $B = \mbox{supp}h.$
We then note that $\| B|x|\| < 1$ by uniform convexity.
Also
\begin{eqnarray*}
 |x^*|(B|x|) & = & \int |x^*| B|x| d\mu \\
             & = & \int_B |x^*| |x| d\mu \\
             & = & \int |h| d\mu \\
             & = & 1, \mbox{ a contradiction.} \\
\end{eqnarray*}
Also $\mbox{supp}x^*_i = B$ since $X^*$ is uniformly convex.
Now as in \cite{bib:G} we observe that there exists a measurable function $\theta$
of modulus one so that $ x_2^* = \theta x_1^*.$
Indeed define $\theta = \frac{x_2^*}{x_1^*}$ on B and $\theta = 1$ on
$B^c.$
Then
\begin{eqnarray*}
\int |h||\theta| d\mu & = & \int |x_1| |x_2^*| d\mu \\
                      & \leq & \| |x^*_2| \|_{X^*} \| |x_1| \|_X \\
                      & = & 1 \\
\end{eqnarray*}

{\noindent Similarly,}
$\int |h||{\theta}^{-1}| d\mu \leq 1.$
So $$\int |h|\{|\theta| + |{\theta}^{-1}|\}d\mu \leq 2.$$
And since $t + t^{-1} \geq 2$ for $t > 0$ we get
$$\int |h| \{|\theta| + |{\theta}^{-1}|\} d\mu \geq 2 \int |h| d\mu =2.$$
Thus $|\theta| + |{\theta}^{-1}| = 2$, but this cannot happen unless
$|\theta| =1.$
Thus $|x_1^*| = |x_2^*|.$ Now $\mbox{supp}x_i = \mbox{supp}h$
and $h = |x_1^*| \cdot x_1 = |x_2^*| \cdot x_2$ yields that $x_1 = x_2.$
$\bullet$

\vspace{1cm}
We are now ready to give a proof of the main result of this work.

{\noindent \bf Proof of Theorem 2.1:}
Suppose that$X$ contains $\l_{\infty}^n$ uniformly in $n$.
Then $S(X)$ is not homeomorphic to $S(L_1(\osm))$ for any measure
space $\osm.$
Indeed this follows, as in \cite{bib:O.S}, from Enflo's result \cite{bib:E}
 that the sets
\newline $S(l_{\infty}^n), n \in \integers$
cannot be uniformly embedded into $S(L_1).$

For the converse assume that $X$ does not contain  $\l_{\infty}^n$ uniformly
in $n.$ Then $X$ must be order continuous since $X$ does not contain $c_0$
\cite{bib:L.T}.
Then the proof goes as in \cite{bib:O.S}.
 By a theorem of Maurey and Pisier \cite{bib:MP}
X must have a finite cotype $q'.$
Thus $X$ is $q$-concave, in fact for all $q > q'$ (\cite[p 88]{bib:L.T}).
Renorm $X$ by an equivalent norm for which $M_q(X) = 1$ and such
that $X$ has the
 same lattice structure (see \cite[p 54]{bib:L.T}).
Then the 2-convexification $X^{(2)}$ of $X$ in this norm satisfies
$$M_{2q} (X^{(2)}) = 1 = M^2(X^{(2)})$$
(\cite[p 54]{bib:L.T} ).
This implies that $X^{(2)}$ is uniformly convex and uniformly smooth
(\cite[p 80]{bib:L.T}), and so
$$F_{X^{(2)}} : S(L_1(\mu)) \ra S(X^{(2)})$$ is a uniform
homeomorphism by Proposition 2.9.
Therefore $$G_2 \circ F_{X^{(2)}} : S(L_1(\mu)) \ra S(X)$$
is a uniform homeomorphism by Proposition 2.4.
$\bullet$

\vspace{1cm}
\rmk {\bf 2.10}
\cite{bib:O.S} If $S(X)$ is uniformly homeomorphic
to $S(Y)$ then $Ba X$ and $Ba Y$
are uniformly homeomorphic.

\vspace{1cm}
{\noindent \bf Corollary 2.11:}
If $X$ is a separable infinite dimensional Banach lattice then $S(X)$ and
$S(l_1)$ are uniformly homeomorphic if and only if $X$ does not contain
$l_{\infty}^n$
uniformly.

\vspace{1cm}
{\noindent \bf Proof:}
By Theorem 2.1, $S(X)$ is uniformly homeomorphic to
$S(L_1(\mu))$ for some probability space $\osm$
where $L_1(\mu)$ is separable.
By standard representation theorems either $L_1(\mu) \cong l_1$
or $L_1(\mu) \cong (L_1[0,1] \oplus l_1(I))_1$ where $I$ is countable.
So $S(X)$ is uniformly homeomorphic to $S((L_1[0,1] \oplus l_1(I))_1).$
Then one can define
$$H:S((L_1[0,1] \oplus l_1(I))_1) \ra S((l_1 \oplus l_1(I))_1)$$
as follows: Let $F$ be a uniform homeomorphism between $S(L_{1})$ and
$S(l_{1})$. (Such homeomorphism exists by \cite{bib:O.S}).
If $(g,x) \in S(L_1[0,1] \oplus l_1(I))_1$
then define
$H(g,x) = \left(\|g\|F\left(\frac{g}{\|g\|}\right), x \right)$ for
$g \neq 0$
and $H(0,x) = (0,x).$
It is easily checked that $H$ is a uniform homeomorphism
and now, since $I$ is countable,
\newline $l_1 \oplus l_1(I) \equiv l_1$
which proves the Corollary.
$\bullet$

\vspace{1cm}
\rmk {\bf 2.12}:
In \cite{bib:R}, Y.Raynaud already obtained that if the unit ball
of a Banach space $E$, embeds uniformly into a stable Banach space
$F$, then $E$ does not contain $c_0$. He also proved that if $F$ is
supposed superstable then $E$ does not contain $l_{\infty}^n$ uniformly.
Since $L_1$ is superstable, we could get one direction of Theorem 2.1
in the separable case using the result of \cite{bib:R}.

\vspace{1cm}
\rmk {\bf 2.13}:
If $X$ is $q$-concave
with constant $1$, then $X^{(2)}$ satisfies
$$M_{2q}(X^{(2)}) = M^2(X^{(2)}) = 1,$$
(\cite[p 54]{bib:L.T}) and as we noted before, $X^{(2)}$
 is uniformly convex and
uniformly smooth (\cite[p 80]{bib:L.T}).
We then proved that $$F_{X^{(2)}}: S(L_1(\mu)) \ra S(X^{(2)})$$
 is a uniform homeomorphism with modulus of continuity of $F_{X^{(2)}}$
depending only on the modulus of uniform convexity $\delta_{X^{(2)}}(\ep)$
 of $X^{(2)}$
(which in turn is of power type 2, i.e for some constant
\newline$0 < K < \infty , \;\; {\delta}_{X^{(2)}}(\ep) \geq K {\ep}^2.$
 (\cite[p 80]{bib:L.T}))
 and the modulus of continuity
of $(F_{X^{(2)}})^{-1}$ depending only on the modulus of uniform smoothness
${\rho}_{X^{(2)}} (\tau)$ of
$X^{(2)}$ (which in turn is of power $2q$ i.e. for some constant
$0 < K < \infty , \;\;{\rho}_{X^{(2)}} (\tau) \leq K {\tau}^{2q})$
  \cite[p 80]{bib:L.T}.)

\vspace{1cm}
{\noindent \bf Proof of Theorem 2.2:}
We first observe that $X$ and $Y$ must have weak units,
since they are separable \cite[p 9]{bib:L.T}; and are order continuous since
they both don't contain $c_0.$
In fact, since $q < \infty$ and $q' < \infty, X \;\; \mbox{and} \;\; Y$
don't contain $l_{\infty}^n.$ So, by Corollary 2.11, $S(X)$ and $S(Y)$
are uniformly homeomorphic to $S(L_1).$
Let $\bar{X}$ be $X$ endowed with an equivalent norm and
the same order, for which $M_q(\bar{X}) = 1,$
 and let $\bar{Y}$ be $Y$ with an equivalent norm
 and the same order, for which $M_{q'}(\bar{Y}) = 1.$
With the previous notations used throughout this work,
we have the following diagram:
$$S(X) \stackrel{u^{-1}}{\ra} S(\bar{X})
\stackrel{(G_{\bar{X},2})^{-1}}{\ra} S(\bar{X}^{(2)})
\stackrel{(F_{\bar{X}^{(2)}})^{-1}}{\ra} S(L_1)
\stackrel{F_{\bar{Y}^{(2)}}}{\ra} S(\bar{Y}^{(2)})
\stackrel{G_{\bar{Y},2}}{\ra} S(\bar{Y})
\stackrel{v}{\ra} S(Y)$$
where $v$ is a uniform homeomorphism from $S(\bar{Y})$ to $S(Y)$
with
a modulus of continuity $a$ depending solely on $M_{q'}(Y),$
and $u^{-1}$ is a uniform homeomorphism from $S(X)$
to $S(\bar{X})$ with a modulus of continuity $f$ depending only
on $M_{q}(X).$

{\noindent Let}$$F = v \circ G_{\bar{Y},2} \circ F_{{\bar{Y}}^{(2)}} \circ
(F_{{\bar{X}}^{(2)}})^{-1} \circ (G_{\bar{X},2})^{-1} \circ
u^{-1},$$
then $F$ is clearly a homeomorphism and
$$F^{-1} = u \circ G_{\bar{X},2} \circ F_{{\bar{X}}^{(2)}} \circ
(F_{{\bar{Y}}^{(2)}})^{-1} \circ (G_{\bar{Y},2})^{-1} \circ v^{-1}.$$

{\noindent Let} $b,c,d \; \; \mbox{and}\; \; e$ be
respectively the modulus of continuity of
respectively $G_{\bar{Y},2}, F_{{\bar{Y}}^{(2)}}, (F_{{\bar{X}}^{(2)}})^{-1},
(G_{\bar{X},2})^{-1}.$
$b$ and $e$ are functions solely of 2 by Proposition 2.4
while $c$ and $d$ are functions of $q'$ and $q$ by Proposition 2.9,
Proposition 2.8, and the remark 2.13 above.
Then the modulus of uniform continuity $\alpha$ of $F$
is of the form $\alpha = a \circ b \circ c \circ d \circ e \circ f$
and is a function solely of $q, q', M_q(X), M_{q'}(Y).$
Note that the modulus of continuity of $F^{-1}$ is also given
by $a \circ b \circ c \circ d \circ e \circ f.$
$\bullet$

\vspace{1cm}
{\noindent \bf Proof of Theorem 2.3:}
The proof is exactly the same as in Theorem 2.2 with
the only difference that $F = F_Y \circ (F_X)^{-1}.$
Indeed we have now the diagram:
$S(X) \stackrel{(F_X)^{-1}}{\ra} S(L_1) \stackrel{F_Y}{\ra} S(Y).$
We then let $F = F_Y \circ (F_X)^{-1}$ and use Proposition 2.9
to get that the modulus of continuity of $F$ depends solely
on the modulus of uniform convexity of $Y$ and the modulus of
uniform smoothness of $X.$
$\bullet$

\bigskip

F.  Chaatit

Department of Mathematics

The University of Texas at Austin

Austin, TX 78712-1082 U.S.A.

chaatit@math.utexas.edu

\end{document}